\documentclass{article}
\usepackage[a4paper,tmargin=1.5in,bmargin=1.5in]{geometry}

\usepackage{amsmath} 
\usepackage{amssymb}  
\usepackage{color}
\usepackage{booktabs}
\usepackage{mathtools}
\usepackage{subfig}
\usepackage{hyperref}
\usepackage{breakurl}
\usepackage[normalem]{ulem}
\usepackage{caption}
\usepackage[detect-all]{siunitx}
\usepackage{etoolbox}
\usepackage{multirow}
\usepackage{tabulary}
\usepackage{tabularx}

\newtheorem{rema}{Remark}

\newcommand{\R}{\mathbb R}
\newcommand{\beq}{\begin{equation}}
\newcommand{\eeq}{\end{equation}}
\newcommand{\beqs}{\begin{equation*}}
\newcommand{\eeqs}{\end{equation*}}

\renewcommand\Re{\operatorname{Re}}
\newcommand{\imagunit}{\bf{i}}

\newcommand{\matlab}{{MATLAB}}

\newcommand{\hifoo}{{\sc hifoo}}

\newcommand{\hanso}{{\sc hanso}}
\newcommand{\shifoo}{{\sc hifoos}}
\newcommand{\Hifoo}{{\sc Hifoo}} 

\newcommand{\granso}{{\sc granso}}

\newcommand{\hinfstruct}{{\texttt{hinfstruct}}}

\newcommand{\getPeakGain}{\texttt{getPeakGain}}
\newcommand{\eig}{\texttt{eig}}
\newcommand{\eigs}{\texttt{eigs}}

\newcommand{\complib}{$COMPl_eib$}

\newcommand{\hsym}{H}
\newcommand{\lsym}{L}

\newcommand{\hinf}{\hsym_\infty}
\newcommand{\linf}{\lsym_\infty}

\newcommand{\nAone}{n_x}
\newcommand{\nBone}{n_w}
\newcommand{\nBtwo}{n_u}
\newcommand{\nCone}{n_z}
\newcommand{\nCtwo}{n_y}

\newcommand{\Korder}{n_K}

\newcommand{\consym}{K}
\newcommand{\conmat}[1]{#1_\consym}

\newcommand{\Acl}{\mathcal{A}}
\newcommand{\Bcl}{\mathcal{B}}
\newcommand{\Ccl}{\mathcal{C}}
\newcommand{\Dcl}{\mathcal{D}}
\newcommand{\Ared}{\Acl_{\mathrm{r}}}
\newcommand{\Bred}{\Bcl_{\mathrm{r}}}
\newcommand{\Cred}{\Ccl_{\mathrm{r}}}
\newcommand{\Dred}{\Dcl_{\mathrm{r}}}
\newcommand{\Afull}{\Acl_{\mathrm{f}}}

\newcommand{\tf}{G}
\newcommand{\tfHnorm}{\|\tf\|_{\hinf}}
\newcommand{\tfLnorm}{\|\tf\|_{\linf}}
\newcommand{\tfred}{\tf_{\mathrm{r}}}

\newcommand{\tfredLnorm}{\|\tfred\|_{\linf}}
\newcommand{\tffull}{\tf_{\mathrm{f}}}
\newcommand{\tffullLnorm}{\|\tffull\|_{\linf}}

\newcommand{\st}{\mathrm{subject~to}}

\def\toll{\delta_\mathrm{low}}
\def\tolh{\delta_\mathrm{high}}

\makeatletter
\def\blfootnote{\xdef\@thefnmark{}\@footnotetext}
\makeatother

\title{\LARGE \bf
Low-Order Control Design using a Reduced-Order Model with a Stability Constraint on the Full-Order Model\blfootnote{This work was supported in part by the Research Training Group RTG
2297/1, ``Mathematical Complexity Reduction - MathCoRe" in Magdeburg, funded by
Deutsche Forschungsgemeinschaft,  and in part by the U.S. National Science Foundation under grant DMS-1620083
}
}

\author{
Peter Benner\thanks{Peter Benner and Tim Mitchell are with the Max Planck Institute for Dynamics of Complex Technical Systems,
39106 Magdeburg, Germany {\tt\small benner@mpi-magdeburg.mpg.de, mitchell@mpi-magdeburg.mpg.de}}  \and 
Tim Mitchell\footnotemark[1]  \and
Michael L.~Overton\thanks{Michael L. Overton is with the Courant Institute of Mathematical Sciences, New York University, New York, NY 10012, USA
        {\tt\small overton@cims.nyu.edu}}
}
\date{Feb. 28, 2018}

\begin{document}

\maketitle

\begin{abstract}  

We consider low-order controller design for large-scale linear time-invariant dynamical systems with inputs and outputs.
Model order reduction is a popular technique, but controllers designed for reduced-order models
may result in unstable closed-loop plants when applied to the full-order system. 
We introduce a new method
to design a fixed-order controller by minimizing the $\linf$ norm of a reduced-order closed-loop transfer matrix 
function subject to 
stability constraints on the closed-loop systems for both the reduced-order and the full-order models.
Since the optimization objective and the constraints are all nonsmooth and nonconvex we use a sequential quadratic
programming method based on quasi-Newton updating that is intended for this problem class,
available in the open-source software package \granso.
Using a publicly available test set, the controllers obtained by the new method are compared with those computed by the
\hifoo\ (H-Infinity Fixed-Order Optimization) toolbox applied to a reduced-order model alone,
which frequently fail to stabilize the closed-loop system for the associated full-order model.
\end{abstract}

\section{Introduction}

As the size of dynamical systems continues to grow rapidly, reduced-order modeling 
\cite{morAnt05,morBauBF14,morBenMS05,morBenCOW17}
has become essential.  However, straightforward control design using reduced-order models 
may result in unstable closed-loop plants for the full-order model.  Though there exist methods guaranteeing closed-loop stability
for the latter, in particular frequency-weighted balancing techniques that perform a combined plant-and-controller reduction \cite{morObiA01}, these techniques are often challenging from a computational perspective. 
In this paper we consider designing \emph{low-order} controllers via $\hinf$ optimization applied to the closed-loop system for a \emph{reduced}-order model (ROM), subject
to a stability constraint on the closed-loop system for the associated \emph{full}-order model (FOM).

Consider the open-loop linear time invariant dynamical system \cite{ZhoDG96}
\beq \label{eq:openloop}
	\left[
		\begin{array}{c}
		\dot{x} \\ \hline
		z \\ 
		y 
		\end{array}
	\right]
	=
	\left[
		\begin{array}{c|cc}
		A _1 & B_1 & B_2 \\ \hline
		C_1 & D_{11} & D_{12} \\ 
		C_2 & D_{21} & D_{22} 
		\end{array}
	\right]
	\left[
		\begin{array}{c}
		x \\ \hline 
		w \\ 
		u 
		\end{array}
	\right]
\eeq
where $x~\in~\R^{\nAone}$ contains the states, $u~\in~\R^{\nBtwo}$ is the physical (control) input, $w~\in~\R^{\nBone}$ is the performance input, $y~\in~\R^{\nCtwo}$ is the physical (measured) output, $z~\in~\R^{\nCone}$ is the performance output, and the matrices
are real with compatible dimensions. We wish to
design a controller $\consym$ defining
\[
	\left[
		\begin{array}{c}
		\conmat{\dot{x}} \\
		u
		\end{array}
	\right]
	= K
	\left[
		\begin{array}{c}
		\conmat{x} \\
		y
		\end{array}
	\right]
	=
	\left[
		\begin{array}{cc}
		\conmat{A} & \conmat{B} \\
		\conmat{C} & \conmat{D}
		\end{array}
	\right]
	\left[
		\begin{array}{c}
		\conmat{x} \\
		y
		\end{array}
	\right]
\]
where $\conmat{x}~\in~\R^{\Korder}$ is the controller state and $\conmat{n}$ is the order of the controller, 
resulting in the closed-loop system:
\[
	\label{eq:closed_loop}
	\left[
		\begin{array}{c}
		\dot{x}\\
		\conmat{\dot{x}}\\ \hline
		z
		\end{array}
	\right]
	= 
	\left[
		\begin{array}{c|c}
		\Acl & \Bcl \\ \hline
		\Ccl & \Dcl
		\end{array}
	\right]
	\left[
		\begin{array}{c}
		x\\
		\conmat{x}\\ \hline
		w
		\end{array}
	\right] 
\]
with, assuming $D_{22} = 0$ for convenience:
\beq
\label{eq:ABCDdefs}
\begin{aligned}
		\Acl & = 
		\left[
			\begin{array}{cc}
			A_1 + B_2 \conmat{D} C_2 & B_2 \conmat{C} \\
			\conmat{B} C_2 & \conmat{A}
			\end{array}
		\right] \\ 
		\Bcl & = 
		\left[
			\begin{array}{c}
			B_1 + B_2 \conmat{D} D_{21} \\
			\conmat{B} D_{21}
			\end{array}
		\right] \\ 
		\Ccl & = 
		\left[
			\begin{array}{cc}
			C_1 + D_{12} \conmat{D} C_2 & D_{12} \conmat{C}
			\end{array}
		\right] \\ 
		\Dcl & = 
		\left[
			\begin{array}{c}
			D_{11} + D_{12} \conmat{D} D_{21}
			\end{array}
		\right] .
\end{aligned}
\eeq
When $D_{22}$ is nonzero, the formulas for $\Acl,\Bcl,\Ccl,\Dcl$ are not affine, as they involve $(I - \conmat{D}D_{22})^{-1}$
\cite[p.~446]{ZhoDG96}; the condition $\|\conmat{D} D_{22}\|_2 <1$ ensures that this last matrix is well defined.
However, since $D_{22}$ is zero in all problems in
our test set described below, this constraint will not play a role in this paper.

The closed-loop transfer matrix function
\beq   \label{transferfun}
    \tf(s) = \Ccl(sI - \Acl)^{-1} \Bcl + \Dcl
\eeq
maps the performance input $w$ to the performance output $z$.

\section{Low-Order Controller Design}

Let $\alpha:\R^{n\times n}\to \R$ denote the spectral abscissa, defined for a matrix $M$ by
\[
    \alpha(M) = \max\{\Re \lambda: \det(\lambda I - M)=0\}.
\]
A key requirement of a controller $K$ is that the closed-loop system is stable, i.e., $\alpha(\Acl)<0$.
The $\linf$ and $\hinf$ norms of the transfer matrix function are defined by
\beq
	\label{eq:ntf}
    \tfLnorm = \sup_{\omega\geq 0} \|\tf(\imagunit\omega)\|_2
\eeq
and 
\[
   \tfHnorm = \left\{\tfLnorm \mathrm{~if~} \alpha(\Acl)<0;~ \infty \mathrm{~otherwise} \right\}.
\]
Note that the $\linf$ norm is finite as long as $\Acl$ has no eigenvalues on the imaginary axis while the $\hinf$ norm is finite provided that 
the eigenvalues of $\Acl$ are all in the open left half-plane.
The standard method for computing $\tfLnorm$ is the Boyd-Balakrishnan-Bruinsma-Steinbuch (BBBS) algorithm \cite{BoyB90,BruS90},
implemented in the \matlab\ function \getPeakGain.
The BBBS method converges quadratically to a global maximizer of \eqref{eq:ntf} but it involves computing the eigenvalues of a sequence of Hamiltonian matrices of order $2(\nAone+\Korder)$;
the algorithm thus requires $O((\nAone+\Korder)^3)$ operations per iteration. 
Computing the eigenvalues of $\Acl$ then determines $\alpha(\Acl)$
and hence whether $\tfHnorm$ equals $\tfLnorm$ or $\infty$;
checking stability via the \matlab\ function \eig\ is generally at least an order of magnitude faster 
than using the BBBS method to compute $\tfLnorm$.

If $\conmat{n}=\nAone$, methods to compute a controller that minimizes $\tfHnorm$ are well known.
However, it is often desirable to design a \emph{low-order} controller with $\conmat{n}\ll \nAone$;
note that low-order here refers to the size of the controller $K$,
not to whether $\nAone$, the dimension of the state space, has been reduced.
The open-source toolbox \hifoo\ \cite{hifoo,BurHLetal06}, 
dating from 2006 and based in part on \cite{BurHLetal06a},
addresses this problem by employing nonsmooth unconstrained optimization techniques to
minimize $\tfHnorm$ over the controller variable $K$.
Similarly, the \matlab\ code  \hinfstruct\ \cite{hinfstruct} introduced in release R2010b 
and based in part on \cite{ApkN06,ApkN06a},
also optimizes $\tfHnorm$ over the controller variable $K$.
As the \hinfstruct\ code is proprietary, we focus here on \hifoo, 
whose code can be readily modified.
\Hifoo\ employs a two-phase algorithm,
first reducing $\alpha(\Acl)$ by varying the controller $K$ until $\alpha(\Acl)<0$. 
Then,  having obtained a stabilizing controller for which $\tfHnorm$ is finite,
\hifoo\ proceeds to the second phase, locally minimizing $\tfHnorm$, using that stabilizing controller as a starting point.
\Hifoo\ relies on the nonsmooth unconstrained optimization code \hanso\ \cite{hanso} in both phases. 
During the second phase, 
if a trial value of the controller $K$ results in $\alpha(\Acl)$ being nonnegative, the value of $\tfHnorm$ is $\infty$ by
definition, so the line search in \hanso\ rejects this point and decreases the step length.
This strategy guarantees that all iterates accepted by the line search result in finite $\hinf$ norm values.
From its inception, the motivation for \hifoo\ was ease of use, so that an engineer could easily call it to design low-order
controllers without any need to understand details of the optimization method on which it relies. References to many
applications where \hifoo\ has been used appear in \cite{MitO15}.

However, neither \hifoo\ nor \hinfstruct\ is practical for large-scale systems because of the high cost of computing $\tfLnorm$ and hence $\tfHnorm$. 
For this reason, Mitchell and Overton \cite{MitO15} introduced a new experimental code
\shifoo\ that optimizes an \emph{approximation} to $\tfHnorm$  for a large sparse system,
using a recently proposed, scalable algorithm called
hybrid-expansion-contraction (HEC) \cite{MitO16}.
The HEC algorithm is guaranteed to find a lower bound for $\tfHnorm$ and, under reasonable assumptions,
converges to a stationary point (usually at least a local maximizer) of \eqref{eq:ntf}.
In \cite[Table~2]{MitO15}, this new approach consistently led to stable closed-loop systems for the large-scale system, contrasting with the ROM-only-based \hifoo\ approach, which resulted in 5 of 12 controllers that failed to stabilize the original full-order models.
However, approximating $\tfHnorm$ for the large-scale system sometimes led to 
inconsistencies from one controller $K$ to another, 
effectively implying that the function being optimized was discontinuous in $K$.
The result was that optimization would not infrequently halt prematurely due to failure in the line search,
which assumes the function is continuous.
 
\section{The New Formulation}
To eliminate this difficulty, we take a different course in this paper, designing a low-order controller by optimizing $\tfHnorm$ 
for a ROM while ensuring the stability of the closed-loop
system for the FOM.  More specifically, we consider the optimization problem
\begin{subequations}
\label{eq:optprob}
\begin{align}
     \min_{\conmat{A},\conmat{B},\conmat{C},\conmat{D}}  \,\, \| \tfred & \|_{\linf} \quad \st\  \label{obj} \\
     	\alpha(\Ared) &< 0, \label{constabrom} \\
	\alpha(\Afull) &< 0, \label{constabfom}  
\end{align}
\end{subequations}
where $\tfred$ is the transfer function \eqref{transferfun} 
for matrix quadruple $(\Ared,\Bred,\Cred,\Dred)$ 
defined by the matrices of \eqref{eq:ABCDdefs} \emph{built using the ROM matrices of} \eqref{eq:openloop}, and
$\Afull$ is the matrix $\Acl$ also given in \eqref{eq:ABCDdefs} but 
\emph{built using the FOM matrices of} \eqref{eq:openloop}.

The key idea of the new formulation is that while the optimization objective in \eqref{obj} is the ROM transfer function norm
$\tfredLnorm$,
stability of the FOM closed-loop system, specified by inequality constraint \eqref{constabfom}, is a requirement for the design of the controller.
As we discuss below, assessing the stability of the FOM closed-loop system via computing the spectral abscissa of $\Afull$ can be done much
more efficiently than the computation of the FOM transfer function norm $\tffullLnorm$. 
Furthermore, by also specifying stability of the ROM closed-loop system as an explicit inequality constraint, 
given by \eqref{constabrom}, the optimization 
objective may now be chosen as the $\linf$ norm of $\tfred$ instead of the $\hinf$ norm. Note that solving \eqref{eq:optprob} is
mathematically equivalent to minimizing 
\beq
   F(K) = \left \{  \begin{array} {l} \tfredLnorm \mathrm{~if~} \max\{\alpha(\Ared),\alpha(\Afull)\} < 0\\
                                                   \ \  \ \infty \ \ \ \ \  \mathrm{otherwise}. \end{array} \right . \label{eq:Fdef}
\eeq
Although we use the formulation \eqref{eq:optprob} in the algorithmic development below, it is the final computed values of $F(K)$
that we will report in our evaluation of the algorithms. 

The constrained optimization problem given by \eqref{eq:optprob} has  objective
and constraint functions that are all continuous, though none are convex or smooth.
The $\linf$ norm function in \eqref{obj} is locally Lipschitz, but the spectral abscissa function in \eqref{constabrom}--\eqref{constabfom}
is not locally Lipschitz at a matrix $M$ with an eigenvalue $\lambda$
with $\Re\,\lambda = \alpha(M)$ that has multiplicity two or more.
Since standard constrained optimization software packages 
are not intended for cases where the optimization objective or any of its constraints is nonsmooth, we use 
a recently introduced sequential quadratic programming method based on quasi-Newton (BFGS) updating \cite{CurMO17} 
that is intended for this problem class. This method is implemented in the open-source software package
\granso\ (GRadient-based Algorithm for Non-Smooth Optimization) \cite{granso}.  Extensive experimental results on a suite of challenging 
static output feedback control design problems involving multiple plants
were reported in \cite{CurMO17}, exhibiting very good results compared with three other methods. 
In almost all cases, the objective and constraint functions were nonsmooth, and in many cases even non-locally-Lipschitz,
at the approximate solutions computed by \granso.
Compared to approximations found by the other methods, those obtained by \granso\ were often better
both in terms of reduction in the optimization objective and constraint satisfaction, and also in terms of running time.

As its name suggests, in order for \granso\ to be applied to an optimization problem it needs access to the
gradients of the relevant objective and constraint functions.
The philosophy underlying the use of BFGS for nonsmooth optimization \cite{CurMO17,LewO13} is that,
since locally Lipschitz functions such as the $\linf$ norm and semi-algebraic functions such as the spectral abscissa
are differentiable almost everywhere, it makes sense to compute gradients at optimization iterates.
Although the objective and constraints are often \emph{not} differentiable at stationary points, these points
are not normally encountered by the optimization method except in the limit. A remarkable observation concerning 
the use of BFGS on nonsmooth problems
is that the method rarely if ever converges to non-stationary points, as discussed at length in \cite{LewO13}.
Naturally, two gradients evaluated
at two nearby points may be very different, but it is precisely this property that is exploited by the BFGS updating.
The $\linf$ norm is differentiable with gradient coinciding with the gradient of the largest singular value $\gamma$ of the transfer function at the
frequency $\tilde\omega$ maximizing $\|\tf(\imagunit\omega)\|_2$, provided the maximizer is unique and the largest singular value $\gamma$ is simple.
The formulas for its gradient are well known and involve the corresponding right and left singular vectors for $\gamma$.
Likewise the spectral abscissa $\alpha$ is differentiable at a matrix $M$ provided its eigenvalue $\lambda$ with largest real part is unique
or part of a unique complex conjugate pair of eigenvalues and that $\lambda$ is simple.  The formula for the gradient of the spectral abscissa is also well
known and involves the corresponding right and left eigenvectors for the eigenvalue $\lambda$. \Hifoo's gradient calculations are done via forming the matrix derivatives of \eqref{eq:ABCDdefs} but for large-scale systems, this is expensive in terms of storage and computation.
To overcome such inefficient scaling properties, 
\shifoo\ instead obtained the gradients of $\alpha(\Afull)$ with respect to the controller variable $K$ via differentiating inner products 
defined for the matrices of \eqref{eq:ABCDdefs} \cite[Section~3]{MitO15}; 
we adopt the same approach here.

The method used by the optimization code \granso\ allows the use of infeasible points, i.e., values of the controller $K$ for which either
or both of the stability constraints \eqref{constabrom} and \eqref{constabfom} on the ROM and the FOM respectively are violated. 
However, we found that this worked
poorly, because once an infeasible point was generated, moving towards the feasible region typically resulted
in substantial increase in the optimization objective $\tfredLnorm$, often preventing the algorithm, which weighs information
from the objective and constraints together using a penalty function, from finding another feasible point. 
So, instead we used two alternative methods.

\smallskip
{\bf Algorithm 1.} Use \granso\ in unconstrained mode to first:
\begin{enumerate}
\item[] {\bf Stabilize:} by minimizing $\max\{\alpha(\Ared),\alpha(\Afull)\}$ until a feasible point
for the constraints in \eqref{constabrom} and \eqref{constabfom} is found, and then
\item[] {\bf Optimize:} by switching to a second unconstrained phase that directly
minimizes the function $F(K)$ given in \eqref{eq:Fdef} using \granso.
\end{enumerate}\noindent
Termination takes place if a maximum iteration limit is exceeded in either phase, or if, using its default options,
\granso\ determines  that an approximate stationarity condition is satisfied in the optimization phase.
\smallskip

{\bf Algorithm 2.} Repeat the following \emph{in a loop}:
\begin{enumerate}
\item[(A)] use \granso\ in \emph{unconstrained mode} to minimize $\max\{\alpha(\Ared),\alpha(\Afull)\}$, quitting when a feasible point
for the constraints \eqref{constabrom} and \eqref{constabfom} is found, and continuing to:
\item[(B)] use \granso\ in \emph{constrained mode} to solve \eqref{eq:optprob}, quitting when an iterate is generated for which
either \eqref{constabrom} or \eqref{constabfom} is violated, and returning to step (A).
\end{enumerate}
\noindent
Termination takes place if a cumulative iteration limit is exceeded in either (A) or (B), 
or, using its default options, \granso\ determines that an approximate stationarity condition is satisfied in (B).
Since the restabilizations
done in (A) may increase $\tfredLnorm$, Alg.~2 simply keeps track of the best controller encountered for returning
to the user when the computation is finished.
\smallskip

\begin{rema} The stabilize-then-optimize approach of Alg.~1 is effectively
identical to that used in current versions of \hifoo\, except that (a) the new algorithm uses the \granso\ optimization code instead of
the older code \hanso, and (b) existing versions of \hifoo\ do not use the stability constraint \eqref{constabfom} for the FOM.
As in \hifoo, the line search in the second ``optimize" phase of  Alg.~1 ensures stability is maintained by rejecting
any points where the minimization objective $F(K)$ is infinite due to a stability constraint being violated.
\end{rema}

\begin{rema}
Note that Alg.~2 is not just Alg.~1 done in a loop.  The ``optimize" phase of Alg.~1 simply rejects
controllers $K$ that destabilize either the ROM or FOM (since $F(K)$ is infinite for such $K$) while the (B) phase of Alg.~2 takes advantage of the explicit stability 
constraints in attempt to produce better search directions, i.e., 
ones which simultaneously minimize the $\linf$ norm
and maintain stability. 
\end{rema}

\section{Evaluation}

In order to assess the ability of our new methods to find controllers that minimize the $\hinf$ norm of the ROMs
as much as possible while also stabilizing the FOMs,
we applied it to a selection of large-scale 2D heat flow problems from version 1.1 of \complib\ \cite{compleib}.
We chose the same twelve \texttt{HF2Dx} FOM examples used in the evaluation of \shifoo\ 
\cite{MitO15} because \complib\ already provides corresponding medium-scale ROMs 
for these examples. See Table~\ref{table:test_set} for the list of problems chosen and their full and reduced orders.
For all problems, we elected to compute order 10 controllers.

\begin{table}[t]
\small
\centering
\caption{
Test Set Summary
}
\begin{tabular}{l | cc | cccc  }
\toprule
Problem & 
	$\nAone$ (FOM) & $\nAone$ (ROM) & $\nBone$ & $\nBtwo$ &
	$\nCone$ & $\nCtwo$  \\
\midrule
\texttt{HF2D1} & 3796 & 316 & 3798 & 2 & 3796 & 3 \\
\texttt{HF2D2} & 3796 & 316 & 3798 & 2 & 3796 & 3  \\
\texttt{HF2D5}  & 4489 & 289  & 4491 & 2 & 4489 & 4\\
\texttt{HF2D6}  & 2025 & 289  & 2027 & 2 & 2025 & 4 \\
\texttt{HF2D9}  & 3481 & 484 & 3483 & 2 & 3481 & 2\\
\texttt{HF2D\_CD1}  & 3600 & 256  & 3602 & 2 & 3600 & 2\\
\texttt{HF2D\_CD2}  & 3600 & 256 & 3602 & 2 & 3600 & 2\\
\texttt{HF2D\_CD3} & 4096 & 324 & 4098 & 2 & 4096 & 2\\
\texttt{HF2D\_IS1} & 4489 & 361 & 4491 & 2 & 4489 & 4\\
\texttt{HF2D\_IS2} & 4489 & 361& 4491 & 2 & 4489 & 4\\
\texttt{HF2D\_IS3} & 3600 & 256 & 3602 & 2 & 3600 & 2\\
\texttt{HF2D\_IS4} & 3600 & 256 & 3602 & 2 & 3600 & 2\\
\bottomrule
\end{tabular}
\label{table:test_set}
\end{table}

We used  \getPeakGain\ to compute $\tfredLnorm$, setting its tolerance to $\tolh =10^{-14}$  because we have
observed that its default tolerance of $10^{-2}$ often returns insufficiently accurate results, 
and that even $\toll =  10^{-7}$, the tolerance used by \hifoo\ for its $\hinf$ norm calculations, 
can also sometimes lead to numerical problems.  
This issue of numerical accuracy is important since gradients play a key role in
the optimization procedure, and a given tolerance sometimes results in a computed gradient of  $\tfredLnorm$ that is
significantly more inaccurate than the tolerance might suggest, even when the gradient is well defined at the point.

The spectral abscissa of $\Ared$ is computed by a call to \eig, 
while the spectral abscissa of $\Afull$ is computed using \eigs\ 
which accesses the matrix $\Afull$ only via matrix-vector products.
Assuming the original full-order version of $A_1$ is sparse, as is generally the case, 
$\Afull$ can be cheaply applied as a matrix-vector operator.
Thus, by using \eigs, it is relatively efficient to assess the stability of fairly large systems
for which computing the $\linf$ norm is out of reach. 
Although \eigs\ is not absolutely guaranteed to return the eigenvalues of $\Afull$ with largest real part, it is generally quite reliable 
in practice. Indeed, it was fully sufficient, with appropriate parameter choices, for the aforementioned 
experiments validating \shifoo; see \cite{MitO15} for details.
It is necessary to call \eigs\ twice, 
once for $\Afull$ and once for its transpose, to obtain both the right and left eigenvectors corresponding 
to the eigenvalue with largest real part, as these are needed to compute the gradient of the spectral abscissa,
as noted above.

In Tables \ref{table:hinfvals} and \ref{table:times} below, the columns labeled ``Alg.~1" with subheading ``R+F"
and those labeled ``Alg.~2" with subheading ``R+F" refer to Algorithms 1 and 2 specified above. 
For comparative purposes, we also ran versions of these algorithms that omitted the
stability constraint \eqref{constabfom} on the FOM; the columns labeled ``Alg.~1" with subheading ``R only"
and those labeled  ``Alg.~2" with subheading ``R only" refer to these experiments.  Alg.~1 without the FOM stability constraint is effectively the
same as the algorithm used in existing versions of \hifoo\ but, to account for differences in implementations, we also 
ran version 3.5 of \hifoo\ (using version 2.2 of \hanso) to produce controllers using only the ROM data.
For consistency with \granso,
we disabled \hanso's (expensive) ``gradient sampling" phase and 
set parameters \texttt{opts.normtol} and \texttt{opts.evaldist} both 
to $10^{-6}$.  
We ran \hifoo\ twice, once with its default $\toll = 10^{-7}$ tolerance for the BBBS method 
and again with the tighter $\tolh = 10^{-14}$ tolerance we used for \texttt{getPeakGain} 
in our new code.

For each problem, we randomly generated an initial controller so that all methods/variants 
in our comparison would be initialized at the same point. 
Since by default, \hifoo\ also attempts optimization from three automatically randomly-generated 
controllers (determined by the positive integer parameter \texttt{opts.nrand}), 
we disabled these additional starting points via a slight 
modification to the \hifoo\ code so that it would allow $\texttt{opts.nrand}\coloneqq0$.

All numerical experiments were implemented and run using \matlab\ R2017a on a workstation with an Intel Core i7-6700 (4 Cores @ 3.4 GHz) 
and 16 GB memory.

\begin{table}[t]
\small
\caption{Final Values of $F(K)$. }

\setlength{\tabcolsep}{3pt}
\robustify\bfseries
\center
\begin{tabular}{ l | rr | rr | rr } 
\toprule
\multicolumn{7}{c}{Final values of $F(K)$ in eq.~\eqref{eq:Fdef}}\\
\midrule
\multicolumn{1}{c}{} &
	\multicolumn{2}{c}{\hifoo\ v3.5} & \multicolumn{2}{c}{Alg.~1} & \multicolumn{2}{c}{Alg.~2} \\
\cmidrule(lr){2-3}
\cmidrule(lr){4-5}
\cmidrule(lr){6-7}
\multicolumn{1}{l}{Problem} &
	\multicolumn{1}{c}{$\toll$} & \multicolumn{1}{c}{$\tolh$} & 
	\multicolumn{1}{c}{R only} & \multicolumn{1}{c}{R+F} &
   	\multicolumn{1}{c}{R only} & \multicolumn{1}{c}{R+F} \\
\midrule
\texttt{HF2D1}     & \bfseries 6511.1 & 6512.8 & 6519.4 & 6519.4 & 22928.5 & 9718.7 \\ 
\texttt{HF2D2}     & 5609.6 & 5598.8 & \bfseries 5597.0 & 6292.7 & 5600.6 & 9635.9 \\ 
\texttt{HF2D5}     & 17716.1 & 17403.3 & 17317.2 & 39890.5 & 17292.6 & \bfseries 16847.2 \\ 
\texttt{HF2D6}     & 7400.4 & 7406.0 & 7397.7 & 7383.0 & \bfseries 7370.0 & 7392.6 \\ 
\texttt{HF2D9}     & 60.3 & 60.3 & \bfseries 60.3 & 60.3 & 60.3 & 60.3 \\ 
\texttt{HF2D\_CD1} & $\infty$ & $\infty$ & $\infty$ & 47.4 & $\infty$ & \bfseries 4.6 \\ 
\texttt{HF2D\_CD2} & $\infty$ & $\infty$ & $\infty$ & 48.9 & $\infty$ & \bfseries 48.9 \\ 
\texttt{HF2D\_CD3} & $\infty$ & $\infty$ & $\infty$ & 37.5 & $\infty$ & \bfseries 4.9 \\ 
\texttt{HF2D\_IS1} & 46428.8 & 44777.5 & 44247.1 & 1734082.9 & 42779.3 & \bfseries 42771.5 \\ 
\texttt{HF2D\_IS2} & 10342.4 & 10336.5 & 10309.6 & 10327.5 & \bfseries 10206.7 & 10462.9 \\ 
\texttt{HF2D\_IS3} & $\infty$ & $\infty$ & $\infty$ & 259.3 & $\infty$ & \bfseries 8.5 \\ 
\texttt{HF2D\_IS4} & $\infty$ & $\infty$ & $\infty$ & 23.9 & $\infty$ & \bfseries 6.9 \\ 
\midrule
\multicolumn{7}{c}{Relative differences from the best $F(K)$ values (in bold above)}\\
\midrule
\texttt{HF2D1}     & --- & $0.000$ & $0.001$ & $0.001$ & $2.521$ & $0.493$ \\ 
\texttt{HF2D2}     & $0.002$ & $0.000$ & --- & $0.124$ & $0.001$ & $0.722$ \\ 
\texttt{HF2D5}     & $0.052$ & $0.033$ & $0.028$ & $1.368$ & $0.026$ & --- \\ 
\texttt{HF2D6}     & $0.004$ & $0.005$ & $0.004$ & $0.002$ & --- & $0.003$ \\ 
\texttt{HF2D9}     & $0.000$ & $0.000$ & --- & $0.000$ & 0.000 & $0.000$ \\ 
\texttt{HF2D\_CD1} & $\infty$ & $\infty$ & $\infty$ & $9.297$ & $\infty$ & --- \\ 
\texttt{HF2D\_CD2} & $\infty$ & $\infty$ & $\infty$ & $0.000$ & $\infty$ & --- \\ 
\texttt{HF2D\_CD3} & $\infty$ & $\infty$ & $\infty$ & $6.694$ & $\infty$ & --- \\ 
\texttt{HF2D\_IS1} & $0.086$ & $0.047$ & $0.034$ & $39.543$ & $0.000$ & --- \\ 
\texttt{HF2D\_IS2} & $0.013$ & $0.013$ & $0.010$ & $0.012$ & --- & $0.025$ \\ 
\texttt{HF2D\_IS3} & $\infty$ & $\infty$ & $\infty$ & $29.524$ & $\infty$ & --- \\ 
\texttt{HF2D\_IS4} & $\infty$ & $\infty$ & $\infty$ & $2.453$ & $\infty$ & --- \\ 
\bottomrule
\end{tabular}
\label{table:hinfvals}
\end{table}

\begin{table}[t]
\small
\caption{Wall-Clock Running Times.}
\setlength{\tabcolsep}{3pt}
\robustify\bfseries
\center
\begin{tabular}{ l | rr | rr | rr } 
\toprule
\multicolumn{7}{c}{Wall-clock running times (seconds)}\\
\midrule
\multicolumn{1}{c}{} &
	\multicolumn{2}{c}{\hifoo\ v3.5} & \multicolumn{2}{c}{Alg.\ 1} & \multicolumn{2}{c}{Alg.\ 2} \\
\cmidrule(lr){2-3}
\cmidrule(lr){4-5}
\cmidrule(lr){6-7}
\multicolumn{1}{l}{Problem} &
	\multicolumn{1}{c}{$\toll$} & \multicolumn{1}{c}{$\tolh$} & 
	\multicolumn{1}{c}{R only} & \multicolumn{1}{c}{R+F} &
   	\multicolumn{1}{c}{R only} & \multicolumn{1}{c}{R+F} \\
\midrule
\texttt{HF2D1}     & 5866   & 7065   & 7167   & 10556  & 10264  & 11409  \\ 
\texttt{HF2D2}     & 2373   & 5931   & 5546   & 463    & 6205   & 16034  \\ 
\texttt{HF2D5}     & 3675   & 6579   & 6468   & 737    & 3045   & 10896  \\ 
\texttt{HF2D6}     & 4781   & 5755   & 5547   & 6175   & 6458   & 6075   \\ 
\texttt{HF2D9}     & 451    & 431    & 738    & 523    & 989    & 763    \\ 
\texttt{HF2D\_CD1} & 2708   & 3368   & 3387   & 73     & 1130   & 7872   \\ 
\texttt{HF2D\_CD2} & 2741   & 3184   & 3519   & 92     & 4977   & 137    \\ 
\texttt{HF2D\_CD3} & 6040   & 6684   & 7337   & 295    & 7990   & 11752  \\ 
\texttt{HF2D\_IS1} & 7915   & 10838  & 11988  & 164    & 15056  & 14458  \\ 
\texttt{HF2D\_IS2} & 8068   & 10286  & 10203  & 14450  & 17930  & 22717  \\ 
\texttt{HF2D\_IS3} & 1071   & 2151   & 1295   & 114    & 108    & 2023   \\ 
\texttt{HF2D\_IS4} & 1178   & 1653   & 1485   & 185    & 1359   & 5101   \\ 
\midrule
\multicolumn{7}{c}{Running times relative to \hifoo\ v3.5 ($\toll$)}\\
\midrule
\texttt{HF2D1}     & 1      & 1.20   & 1.22   & 1.80   & 1.75   & 1.94   \\ 
\texttt{HF2D2}     & 1      & 2.50   & 2.34   & 0.19   & 2.61   & 6.76   \\ 
\texttt{HF2D5}     & 1      & 1.79   & 1.76   & 0.20   & 0.83   & 2.96   \\ 
\texttt{HF2D6}     & 1      & 1.20   & 1.16   & 1.29   & 1.35   & 1.27   \\ 
\texttt{HF2D9}     & 1      & 0.96   & 1.64   & 1.16   & 2.19   & 1.69   \\ 
\texttt{HF2D\_CD1} & 1      & 1.24   & 1.25   & 0.03   & 0.42   & 2.91   \\ 
\texttt{HF2D\_CD2} & 1      & 1.16   & 1.28   & 0.03   & 1.82   & 0.05   \\ 
\texttt{HF2D\_CD3} & 1      & 1.11   & 1.21   & 0.05   & 1.32   & 1.95   \\ 
\texttt{HF2D\_IS1} & 1      & 1.37   & 1.51   & 0.02   & 1.90   & 1.83   \\ 
\texttt{HF2D\_IS2} & 1      & 1.27   & 1.26   & 1.79   & 2.22   & 2.82   \\ 
\texttt{HF2D\_IS3} & 1      & 2.01   & 1.21   & 0.11   & 0.10   & 1.89   \\ 
\texttt{HF2D\_IS4} & 1      & 1.40   & 1.26   & 0.16   & 1.15   & 4.33   \\ 
\bottomrule
\end{tabular}
\label{table:times}
\end{table}

\begin{figure}
\begin{minipage}[c]{0.5\linewidth}
\includegraphics[scale=0.4]{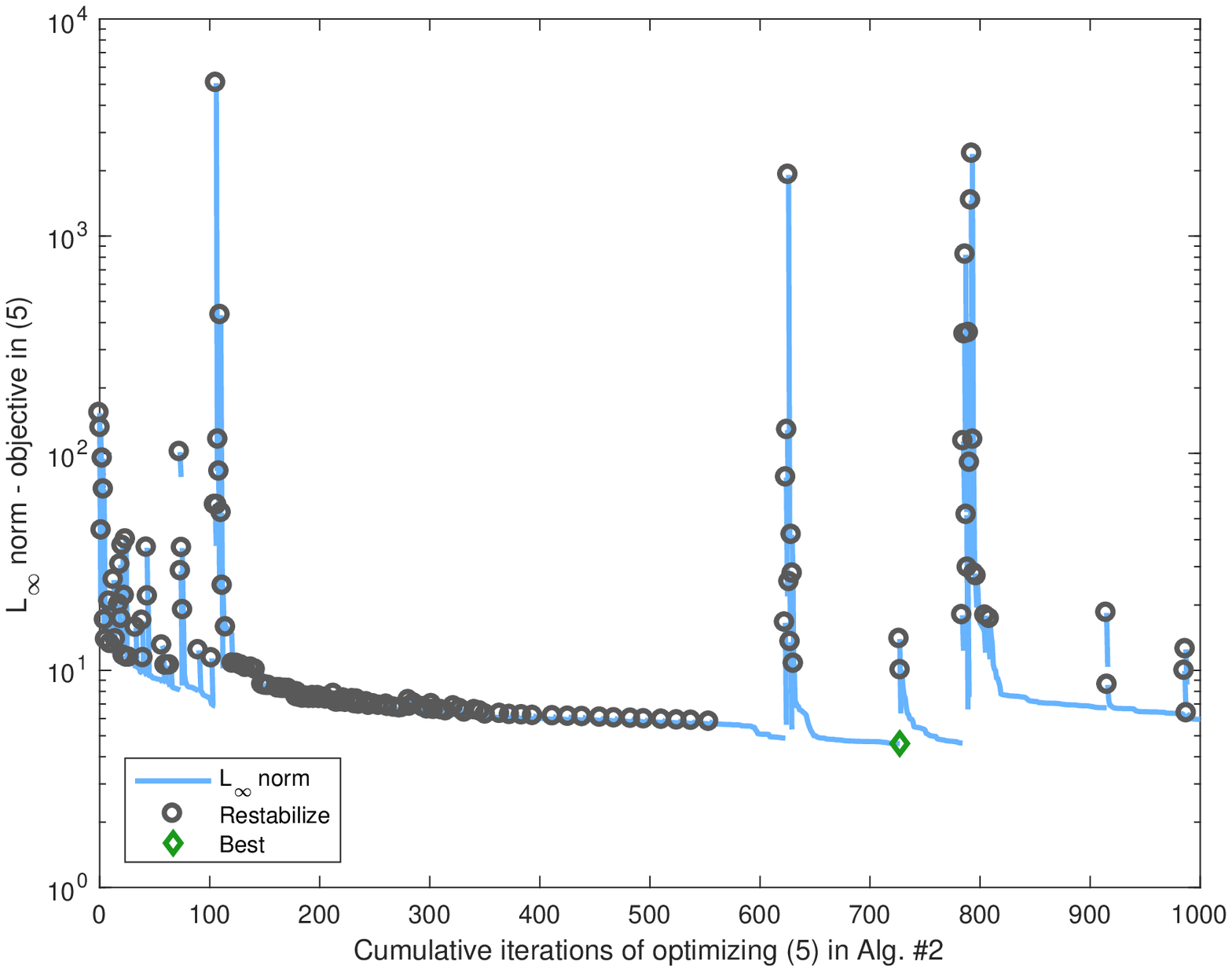}
\caption{Problem \texttt{HF2D\_CD1}}
\label{fig:ex6}
\end{minipage}
\hfill
\begin{minipage}[c]{0.5\linewidth}
\includegraphics[scale=0.4]{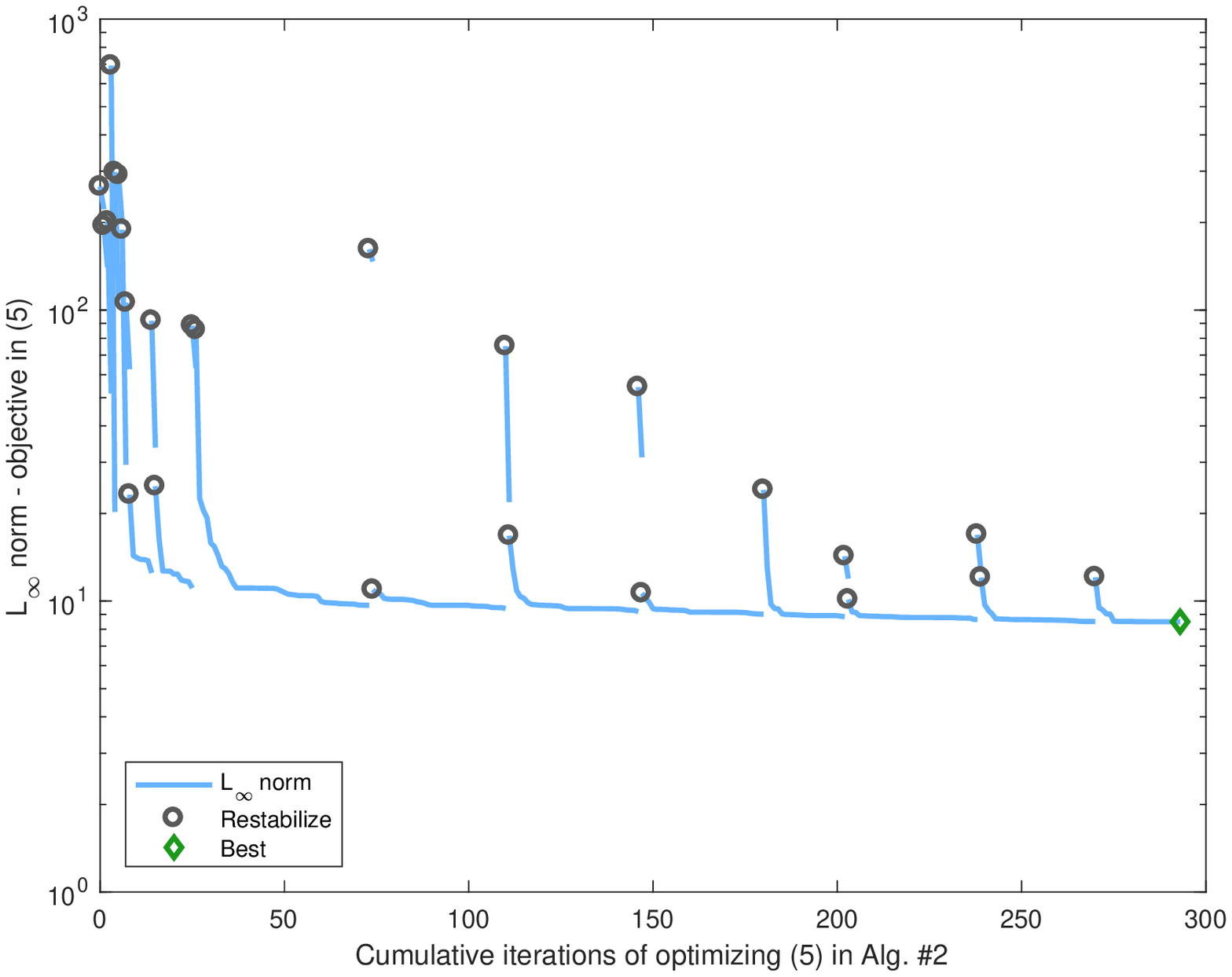}
\caption{Problem \texttt{HF2D\_IS3}}
\label{fig:ex11}
\end{minipage}
\end{figure}

The top half of Table \ref{table:times} reports the total wall-clock running time (in seconds)
for each problem-method pair. The bottom half reports
the ratio of the running times relative to the running times of \hifoo\ v3.5 ($\toll$), with  
values higher than one indicating
how many times slower a method was compared to \hifoo\ while
values less than one indicate the opposite. We see that even though the FOM orders were typically
about 10 times the corresponding ROM orders, the running times for Alg.~2 (``R+F") ranged
from as little as 0.05 to at most 6.76 times the running times of \hifoo\ (which uses only the ROM data),
and less compared to the version of \hifoo\ using the more demanding tolerance $\tolh$ to compute $\tfredLnorm$.

Table \ref{table:hinfvals} reports the best (lowest) values of $F(K)$ (see \eqref{eq:Fdef})
obtained by each method, with, for each problem, the best value obtained over all the methods shown in bold text.
Recall that the value of $F(K)$ is $\infty$ if 
the controller $K$  fails to stabilize \emph{both} the ROM \emph{and} the FOM,
regardless of whether or not the controller $K$ was obtained using any FOM information.
For every ROM-only method in the comparision
(\hifoo\ and the ``R only" variants of Alg.~1 and Alg.~2), we see that $\infty$ is reported in the corresponding columns 
of Table~\ref{table:hinfvals} for 5 out of the twelve problems.
In each case, the respective method's computed controller failed to stabilize the FOM closed-loop systems,
indeed confirming that designing controllers for FOMs, using only ROM information,
can often result in complete failure.  In contrast, both methods that explicitly impose
the FOM stability constraint (the ``R+F" methods) \emph{always} succeeded in simultaneously stabilizing the ROMs and the FOMs,
and hence in Table~\ref{table:hinfvals}, these two ROM-FOM hybrid methods have finite values of $F(K)$ reported for all twelve test problems.  
Furthermore, in addition to ensuring stability of the closed-loop systems for the ROMs and FOMs, Alg.~2 (``R+F")
even succeeded in finding the best value of $F(K)$ on seven of the 12 problems, while on another three,
the values it found were only slightly higher than those obtained by the best methods for those three problems.
This observation is made easier by viewing the bottom half of the table, which shows the relative differences 
of the $F(K)$ values from the best value for each problem,
with dashes indicating that the relevant method was in fact the best.  
A value of $0.000$ means that the relative difference was below our reporting 
limit of $0.001$.

Figs.\ \ref{fig:ex6} and \ref{fig:ex11} show representative examples of the evolution of the
$\tfredLnorm$ values computed by Alg.~2 (``R+F") as a function of the iteration count, for problems
\texttt{HF2D\_CD1} and \texttt{HF2D\_IS3} respectively.
Only the (B) iterations in Alg.~2 are shown as the stabilization iterations in (A) are relatively
less costly. The quantity $\tfredLnorm$ is steadily reduced in (B) until an infeasible point is reached, at which point 
the stabilization phase in (A) typically increases $\tfredLnorm$, sometimes significantly. 
The usual trend, however, is for 
$\tfredLnorm$ to be consistently reduced over a sequence of (A) and (B) iterations, until either:
\begin{itemize}
\item a cumulative maximum number of 1000 iterations in phase (B) is reached, as with problem \texttt{HF2D\_CD1}, or, occasionally, 
\item \granso\ determines that an approximate stationarity condition has been satisfied (see \cite{CurMO17} for details),
as with problem \texttt{HF2D\_IS3}.
\end{itemize}

It is also worth noting that the results for Alg.~1 (``R only") are quite similar to the results for \hifoo, as expected, since these
methods differ only in implementation details. 
Finally, in contrast to the fairly comparable results of Alg.~1 and Alg.~2 in the ``R only" setting, it is clear that Alg.~2 is superior to Alg.~1 in the ``R+F" setting, giving overall very satisfactory results.  
Indeed, the controllers found by Alg.~2 (``R+F") yielded values of $F(K)$ that were on average 8.35 times smaller than those obtained by Alg.~1 (``R+F").

\section{Conclusions}

We have presented a new formulation for minimizing the $L_\infty$ norm of the closed-loop transfer function for a
reduced-order model (ROM) subject to stability
constraints on the closed-loop systems for both the ROM and the full-order model (FOM). 
Algorithm 2 (``R+F") was clearly effective in accomplishing
this goal on the test problems that we considered,
with running times that were not much slower, and sometimes faster,
than the same algorithm without the FOM stability constraint, which, in consequence, often failed to stabilize the 
closed-loop system for the FOM.

\bibliography{csc,mor,software}
\bibliographystyle{alpha}
\end{document}